\documentclass[12pt]{amsart}
\usepackage{latexsym, amsmath,amssymb}

\usepackage{graphicx}

 \numberwithin{equation}{section}

\def\Xint#1{\mathchoice
{\XXint\displaystyle\textstyle{#1}}%
{\XXint\textstyle\scriptstyle{#1}}%
{\XXint\scriptstyle\scriptscriptstyle{#1}}%
{\XXint\scriptscriptstyle%
\scriptscriptstyle{#1}}%
\!\int}
\def\XXint#1#2#3{{\setbox0=\hbox{$#1{#2#3}{%
\int}$ }
\vcenter{\hbox{$#2#3$ }}\kern-.6\wd0}}

\def\dashint{\Xint-}

\renewcommand{\epsilon}{\varepsilon}
\newtheorem{theorem}{Theorem}

\newtheorem{lemma}[theorem]{Lemma}
\newtheorem{corr}[theorem]{Corollary}

\newtheorem{proposition}[theorem]{Proposition}
\newtheorem{deff}[theorem]{Definition}

\newcommand{\bth}{\begin{theorem}}
\newcommand{\ble}{\begin{lemma}}
\newcommand{\bcor}{\begin{corr}}

\newcommand{\bdeff}{\begin{deff}}

\newcommand{\bprop}{\begin{proposition}}
\newcommand{\ele}{\end{lemma}}
\newcommand{\ecor}{\end{corr}}
\newcommand{\edeff}{\end{deff}}

\numberwithin{theorem}{section}

\newcommand{\eprop}{\end{proposition}}

\renewcommand{\Pi}{\varPi}

\renewcommand{\epsilon}{\varepsilon}

\begin{document}

\title[Eigenfunction maxima and spherical means]
{Eigenfunction maxima and spherical means}
%\author[A. C\'ordoba]{Antonio C\'ordoba}
%\address{Instituto de Ciencias Matem\'aticas CSIC-UAM-UC3M-UCM -- Departamento de Matem\'aticas (Universidad Aut\'onoma de Madrid), 28049 Madrid, Spain} 
%\email{antonio.cordoba@uam.es}

\author[A. Mart\'inez]{\'Angel D. Mart\'inez}
\address{Instituto de Ciencias Matem\'aticas (CSIC-UAM-UC3M-UCM) -- Departamento de Matem\'aticas (Universidad Aut\'onoma de Madrid), 28049 Madrid, Spain} 
\email{angel.martinez@icmat.es}

\begin{abstract}
The eigenfunctions of the Laplacian are a central object from the realms of analytic number theory to geometric analysis. We prove that H\"ormander $L^2$-$L^{\infty}$ estimates are equivalent to restriction estimates to small geodesic spheres for a certain class of manifolds.
\end{abstract}

\maketitle

\section{Introduction}

Provided a compact Riemannian manifold $(M,g)$ of dimesion $n$ with Laplace-Beltrami operator $-\Delta_g$ one wishes to understand its eigenfunctions $-\Delta_g\psi_{\lambda}=\lambda^2\psi_{\lambda}$. Recall that due to the elliptic and selfadjoint nature of the operator the eigenfunctions will be smooth and the eigenvalues non-negative real numbers. The following estimate, due originally to L. H\"ormander \cite{H2}, holds:
\[\|\psi_{\lambda}\|_{L^{\infty}(M)}\leq C(M)\lambda^{\frac{n-1}{2}}\|\psi_{\lambda}\|_{L^2(M)}\]
It was extended to similar $L^2$-$L^p$ estimates by C. D. Sogge, whose proof relies on stationary phase estimates for Bessel potentials ocurring in the so-called Hadamard parametrix (see \cite{S1, S2} for details and \cite{D2} for explicit geometrical dependence of the constant). A weaker inequality can be obtained using the Sobolev embedding but the above have the feature of being saturated in the sense that they are sharp when one considers the standard $(n-1)$-sphere $(S^{n-1},\left.g\right|_{\mathbb{R}^n})$.

Our interest on H\"ormander-(Sogge) estimate stems from our work \cite{CM}. The present work originated as an attempt to understand it geometrically without the use of Hadamard's parametrix. Our purpose is to relate the above type of estimate with the following restriction estimate
\[\|\psi_\lambda\|_{L^2(\gamma)}\leq C(M)\|\psi_{\lambda}\|_{L^2(M)}\]
where $\gamma\subset M$ is a geodesic sphere of radius $\approx 1/\lambda$.  During the process we learnt about A. Reznikov work \cite{R1} which was superseded by different methods in the work of Burq-G\'erard-Tzvetkov \cite{BGT}. The former paper considers among other things restriction estimates to geodesic spheres while the latter studies restriction to hypersurfaces with improvements in the case of non-vanishing geodesic  curvature. This encouraged us to better understand the relation between those two problems.

Let us state our result now.

\begin{theorem}\label{Th}%It is equivalent to say that M is a harmonic manifold? See Szabo and Knieper
Let $(M,g)$ be such that the volume of a geodesic ball depends only on its radius and not on its center. Then the following two statements are equivalent:
\begin{enumerate}
\item[(a)] The following estimate holds:
\[\|\psi_{\lambda}\|_{L^{\infty}(M)}\leq C(M,\lambda)\lambda^{\frac{n-1}{2}}\|\psi_{\lambda}\|_{L^2(M)}\]
\item[(b)] There exist a constant $\kappa=\kappa(M)$ such that for any geodesic sphere $\gamma$ of radius $\kappa \lambda^{-1}$ the following estimate holds:
\[\|\psi_{\lambda}\|_{L^2(\gamma)}\leq C(M,\lambda)\|\psi_{\lambda}\|_{L^2(M)}\]
\end{enumerate}
\end{theorem}

The direct implication is obvious, our efforts and presentation will focus in the converse direction. Notice also that the statement itself is quite local in spirit; the radius should be smaller than the injectivity radius to avoid further issues. The scale is {\em the correct one} according to \cite{DF}. One may try to understand the problem globally. As a by-product of H\"ormander estimate one knows $C(M,\lambda)\leq C(M)$ but improvements might be possible. Since it is sharp for the sphere one can not wish better estimates in full generality. However, there has been some results in the negative curvature case; for example, B\'erard proved a logarithmic improvement is possible (cf. \cite{B}). Lot of research has appeared ever since trying to improve this though. %Flat tori are quite interesting since number theory results apply giving $\lambda^{\epsilon}$ estimates for any $\epsilon>0$ in the two-dimensional case. To-date it is not known whether this can be improved to uniform bounds. 
Let us mention that there are results by J. Bourgain \cite{JB} asserting that the exponent of $\lambda$ in H\"ormander's bound can not be improved beyond a certain limit for a perturbation of the two dimensional flat torus. Flat tori are quite interesting since number theory gets into the picture. In the case of constant negative curvature there are considerable improvements on the exponent (cf. Iwaniec and Sarnak \cite{ISa}), open conjectures (cf. Sarnak \cite{Sa}) and counterexamples available in dimension $n=3$ (cf. Rudnick and Sarnak \cite{RuSa}) and $n\geq 5$ (cf. Donelly \cite{D1}). 

Finally, let us mention as a curiosity that A. Zygmund's Theorem 4 \cite{Z} has the same flavour of those conjectures from the restriction side of our equivalence. Furthermore, notice that on the one hand a small geodesic circle is {\em as curved as it can be} in a hyperbolic space and even {\em more curved} than its euclidean counterpart; on the other hand the current philosophy on Fourier restriction estimates is that curvature lurks behind it: this might not be a coincidence.

\section{A comparison principle}\label{comparison}

This section is a technicality, a minor elaboration of some results from \cite{PW} to our present interests. The reader willing to read the core of the proof may proceed to section \ref{proof} coming back here whenever they need to. We will consider the following second order differential operator
\[Lu=u''(x)+g(x)u'(x)+h(x)u(x)\]
acting on sufficiently smooth functions $u$ defined on a fixed open interval; without loss of generality we will assume it to be $(0,1)$. The functions $g$ and $h$ are supposed to be defined and bounded in any closed subinterval of $(0,1)$.  Maximum principles, as the celebrated Hopf's lemma \cite{H}, can be used to prove uniqueness results, our present interest relies on the following related result:

\begin{proposition}
Let $u$ be such that $Lu\geq 0$ in the interval and suppose the existence of some $\varphi>0$ satisfying $L\varphi\leq 0$. Then $\frac{u}{\varphi}$ can not achieve a positive maximum.\\%If $u$ is non negative a priori, it remains non negative in the interval.\\

Furthermore, if $g(x)$ is positive with a singularity at zero of the type $x^{\alpha}$ for some $\alpha\geq -1$ then $u(x)\geq\varphi(x)$ for any $x\in(0,1)$ provided $u(0)=\varphi(0)$ and $u'(0)=\varphi'(0)$.
\end{proposition}

\textsc{Proof:} the first part can be found in \cite{PW}, we include it here for the sake of completeness. We will not prove it in the most straightforward way so as to stress where the hypothesis are needed and how the proof elaborates over rather simple arguments. The second part is an extension suggested by our present needs and new to the best of our knowledge.

If $h$ is non positive there is no need to use $\varphi$. (In fact, $\varphi\equiv 1$ has the desired properties.) Under such restriction let us suppose, arguing by contradiction, that $u$ is not constant and has a positive maximum at $x_0\in (0,1)$; then, the first derivative vanishes at that point $u'(x_0)=0$ and,  finally, as a consequence of $Lu\geq 0$, one gets $u''(x_0)\geq -h(x)u(x_0)\geq 0$ from the hypothesis. This is in contradiction with the fact that it is a maximum since, in such a case, $u''(x_0)\leq 0$ unless $u''(x_0)=0$. To overcome this difficulty and rule out the possibility that $u''(x_0)=0$ one can use a barrier function as follows: first, notice that the above argument holds if one considers the maximum to be achieved at an endpoint $x_0$ in such a way that $u'(x_0)=0$. (This can be understood as a one-dimensional Hopf's lemma.) We will use this view now: instead of considering $u$ we consider $u+v$ on $(0,x_0)$ for some $v$ satisfying $Lv\geq 0$, $v\leq 0 $, $v(x_0)=0$ and $v''(x_0)>0$. This can be achieved taking $v(x)=e^{-M(x-x_0)}-1$ for an appropiate choice of the constant $M$. Now the one-dimensional Hopf's lemma applied to $u+v$ provides the desired contradiction. (Notice $u+v$ has a maximum at $x_0$ if $u$ does.)

Let us recall now that in our statement we had no positivity assumption on $h$. The way to overcome this reduces, precisely, to the existence of $\varphi$ as in the statement. Indeed, if one considers $u=v\varphi$ it is easily proved that 
\[0\geq Lu=L(v\varphi)=\varphi v''+(2\varphi'+g\varphi)v'+L\varphi\cdot v\]
which shows such a $v$ (which is well-defined under our hypothesis) satisfies an equation of the same type with $\tilde{g}=2\frac{\varphi'}{\varphi}+g$ and $\tilde{h}=\frac{1}{\varphi}L\varphi\leq 0$ and hence satisfying the special hypothesis for $v$.

Let us now go over the last part of the statement. We will proceed as before presenting first a more transparent version of the proof and elaborating it to the more refined one. Let $w=u-\varphi$ which will satisfy $Lw\geq 0$, $w(0)=0$ and $w'(0)=0$. One might deduce from the differential inequality, using Taylor's expansion and the hypothesis, that $w''(0)\geq 0$. If we consider the case when $w''(0)>0$ and $h$ is non positive we will be done,indeed, the solution would increase initially; as a consequence: $w>0$ in $(0,\epsilon)$ for some small $\epsilon>0$. This would end the argument from the non existence of a positive maximum as follows: if it becomes negative at some stage it should achieve a maximum meanwhile. Contradiction. But $w''(0)$ might vanish and $h$ might be positive. To dispose of this generalities we will define instead $w_{\delta}=u-\varphi+\delta v$ for some small positive $\delta$ and a barrier $v$ satisfying $Lv\geq0$, $v\geq 0$ and $v''>0$. A function with this properties exist, e.g. $v(x)=e^x-x-1$ as can be checked. (This is where the positivity of $g$ is crucial.) Notice that for such a function $w''_{\delta}>0$ holds and, as before, it increases for some $(0,\epsilon(\delta))$. As a by-product it remains positive and the same is true for the quotient $w_{\delta}\varphi^{-1}$; which can not achieve a maximum due to the first part of our result, as a consequence, it remains positive $w_{\delta}\geq 0$ for any $\delta$. Taking $\delta$ tend to zero concludes the argument.

\section{Proof of theorem \ref{Th}} \label{proof}

The proof will reduce to the study of certain spherical means and hence it is related to, but not subsummed in, Hadamard's parametrix method (cf. \cite{S2} and \cite{J}). Let us sketch the argument first: we will study the spherical means of a smooth function $f$, namely:
\[I_f(x,r)=\dashint_{\partial B_r(x)} f(y)d\sigma(y)\]
It satisfies an explicit second order ordinary differential inequality when $f=\psi_{\lambda}^2$ at the point $x=x_{\lambda}$ where it attains its maximum. This permits us to  compare $I_{\lambda}(x_{\lambda},0)=|\psi_{\lambda}(x_{\lambda})|^2=\|\psi_{\lambda}\|^2_{L^{\infty}(M)}$ with  $I(x_{\lambda},\kappa\lambda^{-1})$ for certain fixed quantity $\kappa$. Nothing else is needed to prove the theorem since for small radius (i.e. $\lambda\gg1$) one can compare the riemannian volume with the euclidean one. The crucial step is based on the comparison principle presented in section \ref{comparison} which enables us to bound $I_{\lambda}(x_{\lambda},r)$ below by some function satisfying an ordinary differential equation. The argument is rather involved since the ordinary differential equation that arises and, as a consequence, its solution depend on $\lambda$; one has to get rid of this dependence so as to find uniform bounds, this is done employing a WKB method for the equation at hand allowing to reduce the argument to a fixed Bessel function if $\lambda$ is big enough. Before proceeding to the proof let us remark, leaving details to the reader, that it is enough to prove the estimate for $\gamma$ a geodesic sphere around a point where the maximum of $|\psi_{\lambda}(x)|$ is achieved. 

Given an eigenfunction $\psi_{\lambda}$ one can consider $I_{\lambda}(r)=I(x_{\lambda},r)$ the spherical means of its square, it satisfies the Euler-Poisson-Darboux differential equation (cf. Appendix) which involves
\[\begin{aligned}
\Delta_gI(x,r)&=\frac{1}{h(r)}\Delta_g\int_{\partial B_r(0)}\psi_{\lambda}(\exp_x(y))^2d\sigma(y)\\%cf. Szabo work on Lichnerowicz conjecture, lemma 1.1
&=\dashint_{\partial B_{r}(0)}\Delta_g(\psi_{\lambda}^2)d\sigma(y)\\
&=2\dashint_{\partial B_{r}(0)}(|\nabla_g\psi_{\lambda}|^2-\lambda^2\psi_{\lambda}^2)d\sigma\\
&\geq-\lambda^2I(x,r)\end{aligned}\]
as a consequence the following differential inequality is satisfied:
\[\frac{d^2}{dr^2}I_{\lambda}(r)+g(r)\frac{d}{dr}I_{\lambda}(r)+2\lambda^2I_{\lambda}(r)\geq 0\]
(See the Appendix to learn more about $g$.) Using our comparison result from section \ref{comparison} one can compare such a function with $J_{\lambda}$ a solution of
\[\frac{d^2}{dr^2}J_{\lambda}(r)+g(r)\frac{d}{dr}J_{\lambda}(r)+2\lambda^2J_{\lambda}(r)=0\]
satisfying $J_{\lambda}(0)=I_{\lambda}(0)=\|\psi_{\lambda}\|^2_{L^{\infty}(M)}$ and $J'_{\lambda}(0)=0$. It will be enough to prove the existence of a constant $\kappa=\kappa(M)$, independent of $\lambda$, such that 
\[J_{\lambda}(r)\geq\frac{1}{2}J_{\lambda}(0)\textrm{ for any }r\in(0,\kappa\lambda^{-1})\]
We will make now the change of variables $\rho=r\sqrt{2}\lambda$ and $\epsilon=(\sqrt{2}\lambda)^{-1}$ to express the above in a more convenient form, namely
\[\frac{d^2}{d\rho^2}K_{\lambda}(\rho)+\frac{n-1}{\rho}\frac{d}{d\rho}K_{\lambda}(\rho)+K_{\lambda}(\rho)=\epsilon k(\rho\lambda^{-1})\frac{d}{d\rho}K_{\lambda}(\rho)\]
where $K_{\lambda}(\rho)=J_{\lambda}(\rho(\sqrt{2}\lambda)^{-1})$ and $k$ is a bounded function that vanishes at zero (cf. Appendix). At least formally this admits a solution of the form
\[K_{\lambda}(\rho)=J_{\lambda}(0)\sum_{n=0}^{\infty}\epsilon^nv_n(\rho)\]
where $v_0$ does not depend on $\lambda$ and satisfies
\[\left\{\begin{array}{l}
v_0''(\rho)+\frac{1}{\rho}v'_0(\rho)+v_0(\rho)=0\\
v_0(0)=1\\
v_0'(0)=0
\end{array}\right.\]
and the rest of the expansion follow interatively from $v_0$:
\[\left\{\begin{array}{l}
v_{n+1}''(\rho)+\frac{1}{\rho}v_{n+1}'(\rho)+v_{n+1}(\rho)=k(\epsilon\rho)v'_n(\rho)\\
v_{n+1}(0)=0\\
v_{n+1}'(0)=0
\end{array}\right.\]
we are abusing notation since there is a hidden dependence in $\lambda$ due to its appearance in $k$. This will not be a problem since all the properties we need from such a function can be shown to be uniform in $\lambda$; namely, that it and its first derivative are uniformily bounded. To make this rigorous one needs to show $v_n$ exist and that the series defining $K_{\lambda}(\rho)$ converges appropiately in an interval for some $\epsilon$ small enough. Once this is done it is clear that for large $\lambda$ the parameter $\epsilon$ will be small and hence one may choose $\kappa$ to be so small that $v_0(\rho)\geq\frac{3}{4}$ holds for any $\rho\in(0,\kappa)$ and then $\lambda$ so large that the error term is smaller than, say, $\frac{1}{4}$. After a change of variables this would end the proof.

The last step will be a consequence of some well-known results from the theory of second order ordinary differential equations (cf. Ince and Sneddon \cite{IS}, chapter 5 for further details). One can write a solution to a general second order ordinary differential equation $Lu=f$ with boundary values $u(0)=u'(0)=0$ as
\[u(x)=y_2(x)\int_0^x\frac{y_1(t)f(t)}{W(y_1,y_2,t)}dt-y_1(x)\int_0^x\frac{y_2(t)f(t)}{W(y_1,y_2,t)}dt\]
where $W$ denotes the wronskian of $y_1$, $y_2$ independent solutions of the homogeneous equation. In our case at hand those solutions are related to some Bessel functions, one of them singular at zero. It will be at zero where we will have to be more careful then. Our claim is that the operation $v_n\mapsto v_{n+1}$ is bounded from $L^{\infty}\rightarrow L^{\infty}$. This would be enough for our purposes and requires further understanding of $y_1$, $y_2$ and their wronskian near zero; namely, we will need that $y_1(x)\approx 1$, 
\[y_2(x)\approx\left\{\begin{array}{cl}
\log(x) &\textrm{if $n=2$}\\
x^{2-n} &\textrm{if $n\geq 3$}
\end{array}\right.\]
and, finally, $W(y_1,y_2,x)\approx x^{1-n}$ for $x\approx 0$. This reduces to knowledge of Bessel functions of the first and second kind since a pair of independent solutions is provided by $y_1(x)=x^{-\frac{n-2}{2}}J_{\frac{n-2}{2}}(x)$ and $y_2(x)=x^{-\frac{n-2}{2}}Y_{\frac{n-2}{2}}(x)$. 

The relevant properties of Bessel functions can be found in Watson's treatise \cite{W}, the standard reference for Bessel functions; in particular we refer to: $\S3\cdot1(8)$, $\S3\cdot51(3)$, $\S3\cdot52(3)$ and $\S3\cdot63(1)$, respectively.

The claim follows from the following equality: 
\[\begin{aligned}
v_{n+1}(x)&=y_2(x)\int_0^x\frac{y_1(t)k(\epsilon t)v'_n(t)}{W(y_1,y_2,t)}dt-y_1(x)\int_0^x\frac{y_2(t)k(\epsilon t)v'_n(t)}{W(y_1,y_2,t)}dt\\
&=\epsilon y_1(x)\int_0^x\frac{y_2(t)k'(\epsilon t)v_n(t)}{W(y_1,y_2,x)}dt-\epsilon y_2(x)\int_0^x\frac{y_1(t)k'(\epsilon t)v_n(t)}{W(y_1,y_2,t)}dt\\
&\qquad\qquad\quad+y_1(x)\int_0^x\frac{d}{dt}\left(\frac{y_2(t)}{W(y_1,y_2,t)}\right)k(\epsilon t)v_n(t)dt\\
&\qquad\qquad\qquad\qquad-y_2(x)\int_0^x\frac{d}{dt}\left(\frac{y_1(t)}{W(y_1,y_2,t)}\right)k(\epsilon t)v_n(t)dt
\end{aligned}\]
where we have used an integration by parts to get rid of $v'_n$ and the aforementioned properties of Bessel functions. Once that equality is known it is straightforward to bound it to obtain a uniform bound $\|v_{n+1}\|_{\infty}\leq C\|v_n\|_{\infty}$. (Notice that $\epsilon$ is bounded since it tends to zero).

To end the argument one still needs to check that the series defining $K_{\lambda}(\rho)$ converges in such a way that it is a function of class $C^2$ in $(0,1)$ and that it satisfies the second order ordinary differential equation that originated it. Let us say that this can be done if one could differentiate the series termwise. To justify it one might prove in a similar guise the estimates $\|v_{n+1}'\|_{\infty}\leq C\|v_n'\|_{\infty}$ and $\|v_{n+1}''\|_{\infty}\leq C\|v_n'\|_{\infty}$. We leave the details to the reader.

\section{Further comments}

The same proof works replacing spherical means of $|\psi_{\lambda}|^2$ by different powers, as a consequence another equivalent assertion would be

\begin{enumerate}
\item[(c)] There exist a constant $\kappa=\kappa(M,p)$ such that for any geodesic sphere $\gamma$ of radius $\kappa \lambda^{-1}$ the following estimate holds:
\[\|\psi_{\lambda}\|_{L^p(\gamma)}\leq C(M,\lambda)\lambda^{\frac{(n-1)(p-2)}{2p}}\|\psi_{\lambda}\|_{L^2(M)}\]
\end{enumerate}

which obviously reduces to (b) for $p=2$. The same argument works for the spherical means of $\psi_{\lambda}$ in which case the differential inequality would turn to be an equality. Let us add that an easy consequence of the equivalence is H\"ormander's estimate for the class of compact manifolds under the hypothesis in the statement since the restriction part can be proved. Indeed, one just needs to notice that underlying the estimate there is an operator $T_{\lambda}u=f$ where
\[\left\{\begin{array}{cl}
-\Delta_gu=\lambda^2u&\textrm{in $\Omega$}\\
u=f&\textrm{in $\partial\Omega$}
\end{array}\right. \]
and $\Omega$ is supposed to be a geodesic ball. Then one can prove the estimate $\|f\|_{L^2(\partial \Omega)}\leq\|u\|_{L^2(\Omega)}$, which is enough. Indeed, the adjoint operator $T^*_{\lambda}$ is an extension operator in the sense that $T_{\lambda}T^*_{\lambda}=\frac{1}{2}Id_{L^2(\partial\Omega)}$. As a consequence $T_{\lambda}$ is uniformily bounded. Let us remark that this argument is not sensitive about global properties inherent to $\psi_{\lambda}$.  Alternatively, one might study the adjoint operator by hand using potential theory to conclude its boundedness. We omit further details, leaving them to the interested reader.

\section{Appendix: the Euler-Poisson-Darboux equation}

This section is a straightforward adaptation of F. John's account of the Euler-Poisson-Darboux equation (cf. \cite{J}, pp. 88-89) where it is deduced in the case of the euclidean space. We include it here for the reader's convenience. With the notation of section \ref{proof} and making some abuse of notation (understanding de integration in a fixed chart with normal coordinates):
\[\int_0^rI(x,\rho)h(\rho)d\rho=\int_{B_r(0)}f(\exp_x(y))d\textrm{vol}_g(y)\]
where $h$ is the Radon-Nykodym derivative $\frac{d\textrm{vol}_g}{d\rho}$. In the case of the $n$-euclidean space $h(\rho)=\omega_{n-1}\rho^{n-1}$. Notice $h(0)=0$ in general since its infinitesimal equivalent is $r^{n-1}$. Furthermore it increases initially and is obviously positive. We may apply the Laplace-Beltrami operator in the variable $x$ to the above equation and use of the divergence theorem as follows
\[\begin{aligned}
\int_0^r\Delta_gI(x,\rho)h(\rho)d\rho&=\int_{B_r(0)}\Delta_gf(\exp_x(y))d\textrm{vol}_g(y)\\
&=\int_{\partial B_r(0)}\frac{\partial}{\partial\nu_y}f(\exp_x(y))d\sigma(y)\\
&=h(r)\dashint_{|z|=1}\frac{\partial}{\partial\nu_y}f(\exp_x(r\cdot z))d\sigma(z)\\
&=h(r)\frac{\partial}{\partial r} I(x,r)\end{aligned}\]
Taking derivatives in $r$ one finally gets the Euler-Poisson-Darboux equation
\[\frac{d^2}{dr^2}I(x,r)+g(r)\frac{d}{dr}I(x,r)=\Delta_gI(x,r)\]
where $g(r)$ denotes the logarithmic derivative of $h(r)$ and hence $g(r)-\frac{n-1}{r}$ is a continuous function vanishing at zero.

\section{Acknowledgments}

The author would like to express his gratitude to Antonio C\'ordoba and Charles L. Fefferman for allowing him to discuss the paper with them; he is also thankful to Eric Latorre who pointed out errata in earlier versions. The work was accomplished during a visit to Princeton University, whose hospitality and environment the author is grateful for; the visit was partially supported by Estancia Breve EEBB-I-16-10718 and MTM2011-2281 project of the MCINN (Spain).

\end{document}